\documentclass[reqno,12pt]{amsart}
\usepackage{amscd,amsfonts,mathrsfs,amsthm,amssymb, amsmath,mathtools}
\usepackage{upgreek}
\usepackage{csquotes}
\usepackage{enumerate}
\usepackage{bbm}
\usepackage{multicol}
\usepackage{stmaryrd,color}
\usepackage{array}
\usepackage{ae}
\usepackage{color}
\usepackage{accents}
\usepackage{epsfig}
\usepackage[e]{esvect}
\usepackage[nobysame]{amsrefs}
\usepackage[margin=1in]{geometry}
\usepackage[toc,title,page]{appendix}
\usepackage{hyperref}
\hypersetup{colorlinks=true,linkcolor=blue,filecolor=blue,urlcolor=blue, citecolor=blue}

\let\oldtocsection=\tocsection
\let\oldtocsubsection=\tocsubsection
\let\oldtocsubsubsection=\tocsubsubsection
\renewcommand{\tocsection}[2]{\hspace{0em}\oldtocsection{#1}{#2}}
\renewcommand{\tocsubsection}[2]{\hspace{1em}\oldtocsubsection{#1}{#2}}
\renewcommand{\tocsubsubsection}[2]{\hspace{2em}\oldtocsubsubsection{#1}{#2}}

\def\equationcolor {\color{black}}
\def\textcolor     {\color{black}}

\def\bcoleq    {\begin{equation}\equationcolor}
\def\ecoleq    {\textcolor\end{equation}}
\def\bcoleqn   {\equationcolor\begin{eqnarray}}
\def\ecoleqn   {\end{eqnarray}\textcolor}

\def\S        {\mathbb{S}}

\def\gind{\operatorname{g}}

\def\R{\mathbb{R}}

\newcommand{\disp}{\displaystyle}
\newcommand{\di}{\mathrm{d}}
\newcommand{\eps}{\varepsilon}
\newcommand{\pa}[1]{{\left(#1\right)}}

\newcommand{\set}[1]{{\left\{#1\right\}}}  
\newcommand{\abs}[1]{{\left|#1\right|}}

\DeclareMathOperator*{\Ric}{Ric}

\newtheorem{theorem}{Theorem}[section]

\newtheorem{lemma}[theorem]{Lemma}

\theoremstyle{definition}
\newtheorem*{assumption*}{$\lambda_{1}$-Condition}

\newtheorem{remark}[theorem]{Remark}

\def\pproof#1{\@ifnextchar[\opargproof
{\opargproof[\it Proof of #1.]}}
\def\opargproof[#1]{\par\noindent {\bf #1 }}

\makeatother

\numberwithin{equation}{section}

\begin{document}

\title[Pinching theorems for complete submanifolds in $\mathbb{S}^{n+p}$]{Sharp pinching theorems for complete submanifolds in the sphere}

\author[M. Magliaro]{\textsc{M. Magliaro}}
\author[L. Mari ]{\textsc{L. Mari}}
\author[F. Roing]{\textsc{F. Roing}}
\author[A. Savas-Halilaj]{\textsc{A. Savas-Halilaj}}

\address{Marco Magliaro \newline
Dipartimento di Scienza e Alta Tecnologia,
Universit\`a degli Studi dell' Insubria,
22100 Como, Italy\newline
{\sl e-mail:} {\bf marco.magliaro@uninsubria.it}
}

\address{Luciano Mari \newline
Dipartimento di Matematica ``Federigo Enriques",
Universit\`a degli Studi di Milano,
20133 Milano, Italy\newline
{\sl e-mail:} {\bf luciano.mari@unimi.it}
}

\address{Fernanda Roing \newline
Dipartimento di Matematica ``Giuseppe Peano",
Universit\`a degli Studi di Torino,
10123 Torino, Italy\newline
{\sl e-mail:} {\bf fernanda.roing@unito.it}
}

\address{Andreas Savas-Halilaj\newline
{Department of Mathematics,
Section of Algebra \!\&\! Geometry, \!\!
University of Ioannina,
45110 Ioannina, Greece} \newline
{\sl E-mail:} {\bf ansavas@uoi.gr}
}

\renewcommand{\subjclassname}{  \textup{2000} Mathematics Subject Classification}
\thanks{A. Savas-Halilaj is supported by (HFRI) Grant No:14758. L. Mari is supported by the PRIN project no. 20225J97H5 ``Differential-geometric aspects of manifolds via Global Analysis"}
\parindent = 0 mm
\hfuzz     = 6 pt
\parskip   = 3 mm
\date{}

\begin{abstract}
We prove that every complete, minimally immersed submanifold
$f\colon\! M^n \to \mathbb{S}^{n+p}$ whose second fundamental form satisfies $|A|^2 \le np/(2p-1)$, is either totally geodesic, or (a covering of) a Clifford torus or a Veronese surface in $\mathbb{S}^4$, thereby extending the well-known results by Simons, Lawson and Chern, do Carmo \&  Kobayashi from compact to complete $M^n$. We also obtain the corresponding result for complete hypersurfaces with nonvanishing constant mean curvature, due to Alencar \& do Carmo in the compact case, under the optimal bound on the umbilicity tensor. In dimension $n \le 6$, a pinching theorem for complete higher-codimensional submanifolds with non-vanishing parallel mean curvature is proved, partly generalizing previous work of Santos. Our approach is inspired by the conformal method of Fischer-Colbrie, Shen \& Ye and Catino, Mastrolia \& Roncoroni.
\end{abstract}

\maketitle
\setcounter{tocdepth}{1}
\section{Introduction}
Throughout the present work, $\S^{n+p}$ denotes the $(n+p)$-dimensional unit sphere, $\S^{n+p}(r)$ that of radius $r$, and $\S^{n+p}_c$ that of sectional curvature $c$.\par
Let $f:M^n\rightarrow\S^{n+p}$, $n\ge 2$, be an immersed submanifold without boundary. According to a seminal  result due to Simons \cite{simons},
if $M^n$ is compact, minimal, and its second fundamental form $A$ satisfies 
\[
|A|^2\le \frac{np}{2p-1}, 
\]
then either $|A|\equiv 0$ or $|A|^2\equiv \frac{np}{2p-1}$. If $|A|\equiv 0$, then $M^n$ is a great sphere,  and if $|A|^2\equiv \frac{np}{2p-1}$, a characterization was given by Lawson \cite{lawson} (in codimension $p=1$) and by Chern, do Carmo \& Kobayashi \cite{chern}: either $M^n$ is a (Riemannian) covering of the Clifford torus given by the natural embedding
\[
T^{n,k} = \mathbb{S}^k\big(\sqrt{{k}/{n}}\big) \times
\mathbb{S}^{n-k}\big(\sqrt{{(n-k)}/{n}}\big) \hookrightarrow \mathbb{S}^{n+1}, \qquad k \in \{1,\ldots, n-1\},
\] 
or a covering of a Veronese surface in $\mathbb{S}^4$. The latter is the embedding of the projective plane $\mathbb{S}^2(\sqrt{3})/\mathbb{Z}_2 \to \mathbb{S}^4$ induced by the restriction to $\mathbb{S}^2(\sqrt{3}) \subset \R^3$ of the map 
\[ 
F: \R^3 \to \R^5, \qquad F(x,y,z) = \left( \frac{yz}{\sqrt{3}},\frac{zx}{\sqrt{3}},\frac{xy}{\sqrt{3}}, \frac{x^2-y^2}{2\sqrt{3}}, \frac{x^2 + y^2 - 2z^2}{6} \right), 
\]
taking the quotient by the antipodal map of $\mathbb{S}^2(\sqrt{3})$. We remark that the characterization of the case $|A|^2\equiv \frac{np}{2p-1}$ is local. For compact submanifolds with nonzero, parallel mean curvature, an analogous problem was addressed by Alencar \& do Carmo \cite{alencar} in codimension $1$, and by Santos \cite{santos} if $p \ge 2$.

The goal of the present paper is to extend the above results to complete, possibly non-compact submanifolds. To introduce our theorems and to explain the main issues to overcome, we first examine the case of hypersurfaces more closely. Let $f:M^n\rightarrow\S^{n+1}$ be a compact immersed hypersurface with (normalized) constant mean curvature $H\ge 0$. For short we call such objects {\em CMC hypersurfaces}.
Denote by $\gind$ the Riemannian metric on $M^n$ and by
$$
\Phi=A-H\gind
$$
the traceless part of the second fundamental form of the hypersurface, and suppose that
\[
|\Phi|^2\le b^2(n,H), 
\]
where $b(n,H)$ is the positive root of the polynomial
\begin{equation}\label{polynomial}
P_{(n,H)}(x)=x^2+\frac{n(n-2)}{\sqrt{n(n-1)}}Hx-n(H^2+1).
\end{equation}
Then, Alencar \& do Carmo \cite{alencar} proved that either $|\Phi|\equiv 0$ and $M^n$ is a sphere or $|\Phi|\equiv b(n,H)$ and $M^n$
covers a minimal Clifford torus or a torus of the form
$$
\S^{n-1}(r)\times\S^{1}(\sqrt{1-r^2})\hookrightarrow\S^{n+1},
$$
of appropriate radius $r\in(0,1)$. This particular example is called $H(r)$-{\em torus}. If $H=0$, notice that the conclusions recover those in \cite{simons,lawson,chern} for $p=1$.

All the proofs of the above-mentioned theorems rely on the strong maximum principle applied to $|\Phi|^2$, which in codimension $1$ satisfies the inequality
\begin{eqnarray}\label{lap1}
\Delta |\Phi |^2\ge -2|\Phi |^ 2 P_{(n, H)}(|\Phi |)+2|\nabla \Phi |^ 2,
 \end{eqnarray}
where $P_{(n, H)}$ is the polynomial given in \eqref{polynomial};
see \cite[page 1226]{alencar}. Indeed, the assumption $|\Phi| \le b(n,H)$ implies that
$P_{(n,H)}(|\Phi|) \le 0$, and therefore $\Delta |\Phi |^2 \ge 0$. Since $M^n$ is compact,
$|\Phi|^2$ must be constant and thus $\nabla \Phi \equiv 0$. The conclusion follows from a careful analysis of hypersurfaces with $\nabla \Phi \equiv 0$.

Seeking to obtain the same results under the weaker assumption that $M^n$ is complete, first observe that a computation due to Leung \cite{leung} shows that the Ricci curvature of $\gind$ satisfies
\[
\Ric \ge -\frac{n-1}{n}P_{(n,H)}(|\Phi|)\gind.
\] 
Consequently, under our assumption it follows that $\Ric \ge 0$ on $(M^n, \gind)$. In dimension $n=2$, Bishop-Gromov's theorem implies that $M^2$ has quadratic area growth and is therefore parabolic, see \cite{grigoryan, huber}. Hence, from \eqref{lap1} we again obtain that $M^2$ is either a sphere or (a covering of) a torus.
However, in higher dimensions $\Ric \ge 0$ is not enough to guarantee the parabolicity of $M^n$, and attempts were made to achieve the goal via the Omori-Yau maximum principle at infinity (see \cite{prs,amr} for a thorough investigation of the principle and its geometric applications). Although there are partial
results, the problem remains still open.
As a matter of fact, one can easily show from \eqref{lap1} that $M^n$ is a sphere whenever $\sup|\Phi|<b(n,H)$; see for instance \cite{yanglian} and \cite{vlachos}. Also, if $|\Phi(x_0)|=b(n,H)$ for some $x_0$, then the strong maximum principle implies $|\Phi| \equiv b(n,H)$ and the classification in \cite{alencar} follows. Consequently, the main difficulty is to characterize the complete CMC hypersurfaces of the unit
sphere with
$$
|\Phi|<b(n,H)\quad\text{and}\quad\sup|\Phi|=b(n,H).
$$
In such a case, finding differential inequalities for which the Omori-Yau principle yields useful information seems a hard task.\par
For this reason, we here pursue a different strategy, inspired by the works of Fischer-Colbrie \cite{colbrie}, Shen \& Ye \cite{shen1} and Catino, Mastrolia \& Roncoroni \cite{catino}. The idea
is to conformally change the metric of $M^n$ by a suitable power of the function
$$
u=b^2(n,H)-|\Phi|^2>0,
$$
to show that $M^n$ is compact, from which it readily follows that $M^ n$ is a sphere. We obtain:

\begin{theorem}\label{THME}
Let $f:M^n\rightarrow\S^{n+1}$ be a complete immersed hypersurface with constant mean curvature
$H\ge 0$. Suppose that the square norm $|\Phi|^2$ of the traceless part of the second fundamental form of $M^n$ satisfies
$$
|\Phi|^2\le b^2(n,H),
$$
where $b(n,H)$ is the positive root of the polynomial \eqref{polynomial} $($in particular, $b(n,0) = \sqrt{n}$$)$.

Then, either $|\Phi| \equiv 0$
$($and $M^n$ is a totally umbilic sphere$)$ or $|\Phi| \equiv b(n,H)$. Furthermore, $|\Phi| \equiv b(n,H)$ if and only if:
\begin{itemize}
\item[(a)] $H = 0$ and $M^n$ covers a Clifford torus $T^{n,k}$ for some $k \in \{1,\ldots, n-1\}$;
\smallskip
\item[(b)] $H > 0$, $n \ge 3$ and $M^n$ covers an $H(r)$-torus with $r^2 < (n-1)/n$;
\smallskip
\item[(c)] $H > 0$, $n=2$ and $M^n$ covers an $H(r)$-torus with $r^2 \neq (n-1)/n$.
\end{itemize}
\end{theorem}

\begin{remark}\label{rem_mainthm}
Let us make some comments.
\begin{enumerate}[\rm(1)]
\item
In Theorem \ref{THME} we are implicitly assuming that $M^n$ is $2$-sided if $H \neq 0$. On the other hand, if $H = 0$, then $M^ n$ is not assumed to be $2$-sided.
\smallskip
\item As discussed in \cite{alencar}, according to the chosen orientation, the mean curvature of an $H(r)$-torus is given by either
\[
H = \frac{(n-1) - nr^2}{nr\sqrt{1-r^2}} \qquad \text{or} \qquad H = \frac{nr^2-(n-1)}{nr\sqrt{1-r^2}}.
\]
The choice leading to positive $H$ is the first one if $r^2 < (n-1)/n$, and by direct computation these $H(r)$-tori satisfy $|\Phi| \equiv b(n,H)$. On the other hand, if $r^2 > (n-1)/n$ the choice is the second one, but for $n \ge 3$ a computation gives $|\Phi| > b(n,H)$. Hence, tori with such $r$ do not satisfy the assumptions in our theorem. The different behaviour is due to the linear term in $P_{(n,H)}$ and does not occur if $n=2$, motivating the presence of $H(r)$-tori with any $r^2 \neq (n-1)/n$ in the classification.
\smallskip
\item
To our knowledge, the use of conformal deformations to get compactness/rigidity properties as outlined above first appeared in works by Schoen and Yau \cite{SchoenY}, Fischer-Colbrie \cite{colbrie}, and Shen \& Ye \cite{shen1}. The method was also exploited to investigate CMC hypersurfaces by Shen \& Zhu \cite{shen2}, Cheng \cite{cheng} and Elbert, Nelli \& Rosenberg \cite{elbert}. It is worth-mentioning that in all these latter results it is required that the manifold $M^n$ has dimension $3$, $4$ or $5$.
The absence of a dimensional constraint in our theorem was somehow unexpected to us.
\smallskip 
\item 
In the unpublished \cite{shenye2}, Shen \& Ye obtained general Bonnet-Myers type results based on the conformal method, see also the recent improvement \cite{catino_ronco} by Catino \& Roncoroni. Once the crucial estimates have been obtained, the final part of our argument can be seen within their framework. However, some details depart from those in \cite{shenye2,catino_ronco}, so we provide full proofs.   

\smallskip
\item
The constant $b(n,H)$ in the theorem is sharp. For example, we infer that the complete minimal hypersurfaces constructed by Otsuki
\cite{Otsuki} and do Carmo \& Dajczer \cite{carmo} have bounded squared norm of the second fundamental form
lying in a range containing $n$; namely
$$
\frac{n(n-1)a_0^2}{1-a_0^2}\le |A|^2\le\frac{n(n-1)a_1^2}{1-a_1^2},
$$
where $a_0$ and $a_1$ are constants such that
$$
0<a_0<\frac{1}{\sqrt{n}}<a_1<1;
$$
for more details see \cite[Section 4 \& Remark 2, p.162]{Otsuki}. Similar results for the norm of the
second fundamental form of a CMC hypersurface with two distinct principal curvatures can be found
in \cite{alias2,perdomo}.
\smallskip
\item
For examples of complete CMC hypersurfaces
in the sphere $\S^{n+1}$ obtained by gluing techniques, we refer to the paper of Butscher \cite{butscher}.
\smallskip
\item
There are several sphere-type theorems in the literature for complete submanifolds in space forms
under various pinching conditions on the second fundamental form;
see for example \cite{asperti,dajczer,hasanis,xu1,xu2,xu3,xu4,xu5}. With the notable exception of 
\cite{hasanis}, the compactness either is assumed or it directly follows from the assumptions,
the Gauss equation, and the Bonnet-Myers theorem.
\end{enumerate}
\end{remark}

Next, we move to higher-codimensional submanifolds $f:M^n\rightarrow\S^{n+p}$. In the minimal case, we are able to obtain a neat extension of the result in \cite{simons,lawson,chern} for $M$ complete:
\begin{theorem}\label{teo_minalta}
Let $p\ge 1$ and let $f:M^n\rightarrow\S^{n+p}$ be a complete minimal immersion. If the norm of the second fundamental form $A$ of $M$ satisfies 
\begin{equation}\label{upperA}
\abs{A}^2\le\frac{np}{2p-1},
\end{equation}
then either $\abs{A}\equiv0$ and $M$ is a totally geodesic sphere, or $\abs A\equiv \frac{np}{2p-1}$. In this latter case, one of the following occurs:
\begin{itemize}
\item[$i)$] $p=1$ and $M^n$ covers a minimal Clifford torus $T^{n,k}$ for some $k \in \{1,\ldots, n-1\}$;
\item[$ii)$] $n=p=2$ and $M^2$ covers a Veronese surface in $\S^4$.
\end{itemize}
\end{theorem}

We also consider submanifolds with non-zero, parallel mean curvature vector. In the compact case, the optimal pinching theorem is due to Santos \cite{santos} and extending it to complete submanifolds reveals to be particularly subtle. We refer to Section 3 for precise statements and more detailed comments.

\section{Proof of Theorem \ref{THME}}

We premit the following well-known algebraic lemma. We include a proof for the sake of completeness.
\begin{lemma}\label{lem_kato}
Let $\Phi: \R^n \times \R^n \to \R^p$ be a vector valued symmetric bilinear form with components $\Phi^{\alpha}_{ij}$, $1 \le i,j \le n$, $1 \le \alpha \le p$. Assume that $\sum_k \Phi^\alpha_{kk} = 0$ for each $\alpha$. Then, the norm
\[
\abs{\Phi}^2 \doteq \sum_{\alpha,i,j} (\Phi^\alpha_{ij})^2
\]
satisfies 
\[
\abs{\Phi}^2 \ge \disp\frac n{n-1}\sum_\alpha\sum_{j=1}^n\pa{\Phi^\alpha_{1j}}^2 \ge \frac n{n-1}\sum_\alpha\pa{\Phi^\alpha_{11}}^2.
\] 
\end{lemma}
\begin{proof}
By Cauchy-Schwarz inequality, in our assumptions
$$
\pa{\Phi^\alpha_{11}}^2 = - \pa{ \sum_{i=2}^n \Phi^\alpha_{ii} }^2 \le (n-1)\sum_{i=2}^n\pa{\Phi^\alpha_{ii}}^2,
$$
hence
\[
\begin{array}{lcl}
\abs{\Phi}^2&\ge&\disp\sum_\alpha\sum_{i=1}^n\pa{\Phi^\alpha_{ii}}^2+2\sum_\alpha\sum_{j=2}^n\pa{\Phi^\alpha_{1j}}^2\ge\frac n{n-1}\sum_\alpha\pa{\Phi^\alpha_{11}}^2+2\sum_\alpha\sum_{j=2}^n\pa{\Phi^\alpha_{1j}}^2\\
&\ge&\disp\frac n{n-1}\sum_\alpha\sum_{j=1}^n\pa{\Phi^\alpha_{1j}}^2 \ge \frac n{n-1}\sum_\alpha\pa{\Phi^\alpha_{11}}^2, 
\end{array}
\]
as claimed. 
\end{proof}

Let $f:M^n\rightarrow\S^{n+1}$ be a complete immersed hypersurface with $|\Phi |\le b(n, H).$ Let us
denote by $b_- < 0 < b_+ = b(n,H)$ the two roots of $P_{(n, H)}$ given in \eqref{polynomial}, namely:  
\[
b_\pm = \pm \sqrt{\frac{n^2(n-2)^2H^2}{4n(n-1)}+n(H^2+1)}-\frac{n(n-2)H}{2\sqrt{n(n-1)}}.
\]
From $|b_-|> b_+$ we deduce that, for $x \in [0, b_+]$, 
\begin{equation}\label{pol_roots}
P_{(n, H)}(x) = (x-b_+)(x-b_-) \le (x-b_+)(x+b_+) = x^2-b_+^2,
\end{equation}
whence \eqref{lap1} can be written in the form
\begin{eqnarray}\label{lap2}
\Delta |\Phi |^2\ge 2(b^2-|\Phi|^2)|\Phi |^ 2 + 2|\nabla \Phi |^ 2\ge 0,
 \end{eqnarray}
where, hereafter, $b = b_+=b(n,H)$. As a consequence, the function
\[
u \doteq b^2 - |\Phi|^2, 
\]
satisfies
\begin{equation}\label{eq_u}
u \ge 0 \quad\text{and}\quad \Delta u \le -2 |\Phi|^2 u \qquad \text{on } \, M.
\end{equation}

We distinguish two cases:

{\bf Case 1:} If $u(x_0) =0$ for some $x_0 \in M^ n$, by the strong maximum principle $u \equiv 0$, whence
$$
|\Phi |^2\equiv b^2\quad \text{and} \quad|\nabla \Phi|\equiv 0.
$$ 
The conclusion follows from \cite[p. 1227-1228]{alencar} or \cite[Theorem 4]{lawson}. More precisely, it is shown that the hypersurface has two distinct principal 
curvatures, 
both constants, and that every $x\in M^n$ has a neighbourhood $U$ for which $f(U)$ is a piece of either a Clifford torus $T^{n,k}$ or an $H(r)$-torus (with $r$ in the range stated in the theorem), according to the value of $H$. We show that $f(M^n)$ is a single such 
torus. Indeed, any Clifford or $H(r)$-torus is the zero set of an appropriate real analytic function on
$\mathbb{S}^ {n+1}$; see for instance \cite[Example 3, p.194]{nomizu}. Fix $x\in M^n$ and $U$ as above and let $\Psi :\mathbb{S}^ {n+1}\to \mathbb{R}$ be 
a real analytic function that vanishes on the torus containing $f(U)$. Being $M^ n$ and $f$ real analytic as  well, 
then $\Psi \circ f$ is real analytic and vanishes on $U$, thus it vanishes on the entire $M^n$. This shows that 
$f(M^n)$ is contained in a single torus $\Sigma^n$. As $f : M^n \to \Sigma^n$ is a local isometry and $M^n$ is complete, Ambrose's Theorem guarantees that $f$ is onto and a Riemannian covering, which 
proves our claim.

{\bf Case 2:} Assume now that $u>0$ on $M^n$. Our goal is to prove that $M^n$ must be a totally umbilic sphere. To reach the goal, inspired by \cite{colbrie,catino} we endow $M^n$ with the metric
$$
\overline{\gind}=u^{2\beta }\gind,
$$
where 
\begin{equation}\label{def_beta}
\beta = \left\{ \begin{array}{ll}
\text{any number in $(0,1)$} & \quad \text{if } \, n=2,3, \\[0.2cm]
\disp \frac{1}{n-2} & \quad \text{if } \, n \ge 4.
\end{array}\right.
\end{equation}
Consider a curve $\gamma : [0,a] \to M$ parametrized by $\gind$-arclength $s$, and denote by $\overline{s}$ the $\overline{\gind}$-arclength of $\gamma$. From 
$$
\partial_{\overline{s}}=u^{-\beta}\partial_s\quad\text{and}\quad d\overline{s}=u^ \beta \di s
$$
the length of $\gamma$ in the metric $\overline{\gind}$ is
\[
\ell_{\overline{\gind}}(\gamma) = \int_0^a u^\beta \di s. \vspace{0.3cm}
\]
We split the proof into three claims: 

\noindent \textbf{Claim 1:} {\em Assume that $\gamma$ is a $\overline{\gind}$-geodesic with non-negative second variation of $\overline{\gind}$-arclength. Then, there exist constants $t_0>1$, $c_0>0$ depending on $n,\beta$ such that 
\begin{equation}\label{int51}
c_0\int_0^{a}u^{\beta}\psi^2 \di s \le -2t_0\int_{0}^{a}u^\beta \psi\psi_{ss}\di s, \qquad \forall \, \psi \in C^2_0([0,a]),
\end{equation}
where $C^2_0([0,a])$ is the set of functions $\psi \in C^2([0,a])$ such that $\psi(0)=\psi(a)=0$.}

\begin{proof}[Proof of Claim 1]
From  the second variation formula, along $\gamma$ we have
\begin{eqnarray*}
\int _0^{\overline{s}(a)} \big\{(n-1)(\varphi_{\overline{s}})^ 2-\varphi^ 2\, \overline{\Ric}(\gamma_{\overline{s}}, \gamma_{\overline{s}})\big\}
 \di\overline{s}
\ge
0, \qquad \forall \varphi \in C^2_0([0,a]),
\end{eqnarray*}
or, equivalently,
\begin{eqnarray}\label{second-variation-length}
\int _0^a \big\{(n-1)(\varphi_s)^2-\varphi^2\,
\overline{\Ric}(\gamma_s,\gamma_s)\big\}u^{-\beta}\di s
\ge 0, \qquad \forall \varphi \in C^2_0([0,a]).
\end{eqnarray}
As it is shown in \cite[Appendix, eqn. (14)]{elbert}, along $\gamma$ the following identity holds:
\[
\overline{{\rm Ric}}(\gamma _s, \gamma _s) ={\rm Ric}(\gamma _s, \gamma _s)
-\beta(n-2)(\ln u)_{ss}-\beta\Delta \ln u .
\]
From \eqref{eq_u}, 
\[
\Delta \ln u \le -2|\Phi|^2 - |\nabla \ln u|^2 \le -2|\Phi|^2 - \{(\ln u)_s\}^2, 
\]
whence
\begin{equation}\label{elbert}
\overline{{\rm Ric}}(\gamma _s, \gamma _s) \ge {\rm Ric}(\gamma _s, \gamma _s) + 2\beta |\Phi|^2
-\beta(n-2)(\ln u)_{ss} + \beta \{(\ln u)_s\}^2.
\end{equation}
Let $\{e_1=\gamma_s,e_2\dots,e_n\}$ be a $\gind$-orthonormal frame along $\gamma$. From the Gauss equation, we have that the components of the Riemann curvature
tensor are given by
$$
R_{ijij}=1-\delta_{ij}+h_{ii}h_{jj}-h_{ij}^2,\quad i,j\in\{1,\dots,n\},
$$
where $h_{ij}$ are the components of the second fundamental form of $M^n$. The identity can equivalently be written in the form
$$
R_{ijij}=1-\delta_{ij}+\big(\Phi_{ii}+H\big)(\Phi_{jj}+H)-(\Phi_{ij}+H\delta_{ij})^2,\quad i,j\in\{1,\dots,n\}.
$$
Using the fact that $\Phi$ is traceless, we get
\[
\Ric(\gamma_s,\gamma_s) = \disp (n-1)-\Phi^2_{11}+(n-2)H\Phi_{11}+(n-1)H^2-\sum_{j=2}^n\Phi^{2}_{1j}. 
\]
By Lemma \ref{lem_kato} we have
\begin{equation}\label{eq_refinedkato}
|\Phi|^2\ge \frac{n}{n-1}\sum_{j=1}^n\Phi_{1j}^2 \ge \frac{n}{n-1} \Phi_{11}^2.
\end{equation}
Fix $\tau \in (0,1]$ and $\eps>0$. Let us estimate the term $H\Phi_{11}$ by means of Young inequality and \eqref{eq_refinedkato} as follows
\begin{equation}\label{eq_key}
\begin{array}{lcl}
H\Phi_{11} & = & H (1-\tau)\Phi_{11} + \tau H \Phi_{11} \\[0.3cm]
& \ge & \disp -H (1-\tau)|\Phi|\sqrt{ \frac{n-1}{n}} - \frac{\tau H^2}{2\eps} - \frac{\tau \eps \Phi_{11}^2}{2}.
\end{array}
\end{equation}
Therefore, 
\begin{equation}\label{r11}
\begin{array}{lcl}
\Ric(\gamma_s,\gamma_s) & \ge & \disp (n-1)- (n-2)\sqrt{ \frac{n-1}{n}}(1-\tau)|\Phi|H + \left(n-1- \frac{(n-2)\tau}{2\eps}\right)H^2 \\[0.4cm]
& & \disp - \left(1 + \frac{(n-2)\tau \eps}{2}\right) \Phi^2_{11} -\sum_{j=2}^n\Phi^{2}_{1j}.
\end{array}
\end{equation}
Since $P_{(n,H)}(b) = 0$ and $|\Phi| \le b$, we have
\[
(n-2)\sqrt{ \frac{n-1}{n}} H|\Phi| \le (n-1)(H^2+1) - \frac{n-1}{n}  |\Phi|^2.
\]
Hence,
\begin{equation}\label{r11_2}
\begin{array}{lcl}
\Ric(\gamma_s,\gamma_s) & \ge & \disp \tau(n-1) + \tau\left(n-1- \frac{n-2}{2\eps} \right)H^2 \\[0.4cm]
& & \disp + (1-\tau)\frac{n-1}{n}  |\Phi|^2 - \left(1 + \frac{(n-2)\tau \eps}{2}\right) \Phi^2_{11} -\sum_{j=2}^n\Phi^{2}_{1j}
\end{array}
\end{equation}
and we conclude that
\begin{equation}\label{estimate1}
\begin{array}{lcl}
\overline{{\rm Ric}}(\gamma _s, \gamma _s) & \ge & \disp \tau(n-1) + \tau\left(n-1- \frac{n-2}{2\eps} \right)H^2 \\[0.4cm]
& & \disp + \left(2\beta + (1-\tau)\frac{n-1}{n}\right)  |\Phi|^2 - \left(1 + \frac{(n-2)\tau \eps}{2}\right) \sum_{j=1}^n\Phi^{2}_{1j} \\[0.5cm]
& & -\beta(n-2)(\ln u)_{ss}+\beta \{(\ln u)_s\}^ 2.
\end{array}
\end{equation}
By \eqref{eq_refinedkato}, 
$$
\begin{array}{l}
\disp \left( 2\beta + (1-\tau)\frac{n-1}{n}\right)|\Phi|^2-\left(1 + \frac{(n-2)\tau\eps}{2}\right)\sum_{j=1}^n\Phi_{1j}^2 \\[0.5cm]
\qquad \ge \disp 
\left(\frac{2n\beta}{n-1} -\tau - \frac{(n-2) \tau\eps}{2}\right)\sum_{j=1}^n\Phi_{1j}^2.
\end{array}
$$
We shall specify $\tau$ and $\eps$ so that 
\[
\left\{\begin{array}{l}
\disp n-1 - \frac{n-2}{2\eps} \ge 0,\\[0.4cm]
\disp \frac{2n\beta}{n-1} -\tau - \frac{(n-2) \tau\eps}{2} \ge 0.
\end{array}\right.
\]
Indeed, it is enough to choose 
\[
\eps = \left\{ \begin{array}{ll}
1 & \quad \text{if } \, n=2, \\[0.2cm]
\disp \frac{n-2}{2(n-1)} & \quad \text{if } \, n \ge 3,
\end{array}\right.
\qquad\text{and}
 \qquad \text{$\tau \ $ small enough.}
\]
For such a choice, setting $c_0 = \tau(n-1)>0$, inequality \eqref{estimate1} yields
\begin{equation}\label{estimate2}
\overline{\Ric}(\gamma _s, \gamma _s)
\ge c_0-
\beta(n-2)(\ln u)_{ss}+\beta \{(\ln u)_s\}^ 2.
\end{equation}
%
%

Replacing \eqref{estimate2}
in \eqref{second-variation-length}, we obtain
\begin{equation}\label{int}
(n-1)\int _0^ a(\varphi_s)^2u^{-\beta}\di s
\ge\int_{0}^{a}\varphi^2u^{-\beta}\Big(c_0-\beta(n-2)(\ln u)_{ss}+\beta\{(\ln u)_s\}^2\Big)\di s.
\end{equation}
Integration by parts gives
\begin{eqnarray*}
-\beta\int_0^ a \varphi^2 (\ln u)_{ss} u^{-\beta}\di s
&\!\!\!=\!\!\!&-\int_0^ a \varphi^2 (\ln u^{\beta})_{ss} u^{-\beta}\di s\\
&\!\!\!=\!\!\!&2\beta\int _0^ a \varphi \varphi_s(\ln u)_{s} u^{-\beta}\di s
-\beta^2\int _0^ a \varphi^ 2 \{(\ln u)_s\}^ 2u^{-\beta}\di s,
\end{eqnarray*}
and plugging into \eqref{int} yields
\begin{eqnarray}\label{int2}
&&\hspace{-38pt}(n-1)\int _0^ a(\varphi_s)^2u^{-\beta}\di s\ge c_0\int_{0}^{a}\varphi^2u^{-\beta}\di s\\
&&\,+2\beta(n-2)\int _0^ a \varphi \varphi_su^{-\beta-1}u_{s}\di s+\beta\big(1-\beta(n-2)\big)\int_{0}^{a}\varphi^ 2 u^{-\beta-2}  (u_s)^ 2\di s.
\nonumber
\end{eqnarray}
We follow ideas in \cite{catino} to treat
the integral inequality \eqref{int2}. Set
\begin{equation}\label{def_varphi}
\varphi=u^{\beta}\psi, \qquad \text{with} \quad \psi \in C^2_0([0,a]).
\end{equation}
Then, \eqref{int2} becomes
\begin{eqnarray}\label{int3}
(n-1)\int_{0}^a(\psi_s)^2u^{\beta}\di s&\!\!\!\ge\!\!\!&c_0\int_{0}^{a}\psi^2u^{\beta}\di s\\
&&+\beta(1-\beta)\int_{0}^{a}u^{\beta-2}(u_s)^2\psi^2\di s
-2\beta\int_{0}^{a}u^{\beta-1}u_s\psi\psi_s\di s.\nonumber
\end{eqnarray}
Define
$$
I \doteq \beta \int_0^{a}u^{\beta-1}u_s\psi\psi_s\di s = \frac{1}{2} \int_0^a (u^\beta)_s (\psi^2)_s.
$$
Integration by parts gives
$$
I= -\frac{1}{2}\int_0^{a}u^{\beta} (\psi^2)_{ss}\di s
=-\int_{0}^{\beta}u^{\beta}(\psi_s)^2\di s
-\int_0^{a}u^{\beta}\psi\psi_{ss}\di s.
$$
For every $t>1$ and every $\varepsilon>0$, from Young's inequality we obtain
\begin{eqnarray*}
2 I&\!\!\!=\!\!\!&2 t I+2(1-t)I\\
&\!\!\!=\!\!\!&-2t\int_0^{a}u^{\beta}(\psi_s)^2\di s-2t\int_0^{a}u^{\beta}\psi\psi_{ss}\di s
+2\beta(1-t)\int_0^{a}u^{\beta-1}u_s\psi\psi_s\di s\\
&\!\!\!\le\!\!\!&-2t\int_0^{a}u^{\beta}(\psi_s)^2\di s-2t\int_{0}^{a}\psi\psi_{ss}u^{\beta}\di s\\
&&+\beta(t-1)\varepsilon\int_0^{a}u^{\beta-2}(u_s)^2\psi^2\di s
+\frac{\beta(t-1)}{\varepsilon}\int_0^{a}u^{\beta}(\psi_s)^2\di s.
\end{eqnarray*}
Choosing
$$
\varepsilon=\frac{1-\beta}{t-1}
$$
we obtain
\begin{eqnarray}\label{estimatecatino}
2I&\!\!\!\le\!\!\!&-2t\int_0^{a}u^{\beta}\psi\psi_{ss}\di s
+\beta(1-\beta)\int_0^{a}\psi^2u^{\beta-2}(u_s)^2\di s\nonumber\\
&&+\frac{\beta(t-1)^2-2t(1-\beta)}{1-\beta}\int_0^{a}u^{\beta}(\psi_s)^2\di s.
\end{eqnarray}
Inserting \eqref{estimatecatino} into \eqref{int3}, we arrive at
\begin{equation}\label{int4}
\int_0^{a}u^{\beta}\big\{c_0\psi^2-p(n,t,\beta)(\psi_s)^2+2t\psi\psi_{ss}\big\}\di s\le 0,
\end{equation}
where
\begin{equation*}
p(n,t,\beta)=\frac{\beta(t-1)^2}{1-\beta}-2(t-1) +(n-3).
\end{equation*}
With our choice of $\beta$, we get
\[
p(n,t_0,\beta) \le 0 \qquad \text{where } \, t_0 = \left\{ \begin{array}{ll}
1 + 2\frac{1-\beta}{\beta} & \text{ if } \, n \in \{2,3\}, \\[0.2cm]
n-2 & \text{ if } \, n \ge 4.
\end{array}\right.
\]
Notice that $t_0>1$, whence it is admissible. Then, \eqref{int4}
becomes \eqref{int51}. 
\end{proof}

\noindent \textbf{Claim 2:} {\em The manifold $M^n$ is compact.}

\begin{proof}[Proof of Claim 2] 
Assume by contradiction that $M^n$ is non-compact. We follow \cite{colbrie} and construct the ``smallest divergent curve" in $(M^n,\overline{\gind})$ issuing from a fixed point. For the sake of completeness, we include a simplified and a bit more general statement:

\begin{lemma}\label{lem_fischercolbrie}
Let $(M^n,\overline{\gind})$ be a non-compact Remannian manifold. Then, for each $o \in M$, there exists a divergent curve $\gamma : [0,T) \to M^n$ issuing from $o$ which is a $\overline{\gind}$-geodesic, minimizes the $\overline{\gind}$-length on any compact subinterval of $[0,T)$ and satisfies the following property:
\[
\text{$\overline{\gind}$ is complete} \qquad \Longleftrightarrow \qquad T = \infty.
\] 
\end{lemma} 

{\em Proof.}
The implication $\Rightarrow$ is obvious, since we know that for a complete Riemannian
metric $\overline{\gind}$ every 
divergent $\overline{\gind}$-geodesic is defined on the entire $[0,\infty)$. To prove the converse, we shall 
construct such a geodesic $\gamma$. Consider an exhaustion of $M^n$ by relatively compact, smooth 
open sets $\Omega_j$ containing $o$. Since $\overline{\Omega}_j$ is a smooth compact manifold with boundary, 
there exists a $\overline{\gind}$-minimizing rectifiable curve $\gamma_j : [0,T_j] \to \overline{\Omega}_j$ 
joining $o$ to $\partial \Omega_j$, which we parametrize by $\overline{\gind}$-arclength. Because $\gamma_j$ is a $\overline{\gind}$-geodesic,
$\gamma_j([0,T_j)) \subset \Omega_j$, and since $\overline{\Omega}_j \subset \Omega_{j+1}$, we have that $\{T_j\}$ is a strictly increasing sequence. Then, 
up to a subsequence,
$$\gamma_j \to \gamma : [0,T) \to M^n,$$
smoothly on compact sets as $j \to \infty$. Since each $\gamma_j$ minimizes 
$\overline{\gind}$-length between any pair of its  points, it follows that $\gamma$
is a $\overline{\gind}$-geodesic that minimizes $\overline{\gind}$-length between any pair of its 
points as well. Moreover, let $\sigma: [0,T_{\sigma}) \to M^n$ be any other divergent curve, 
parametrized by $\overline{\gind}$-arclength. Then, for each natural $j$, let $t_j$ be the first time for which $\sigma(t_j) \in \partial \Omega_j$. By the minimality of $\gamma_j$, and having fixed $S > 0$, we have 
for large enough $j$ that
\[
T_{\sigma} = \ell_{\overline{\gind}}(\sigma) \ge \ell_{\overline{\gind}}(\sigma_{|_{[0,t_j]}}) \ge \ell_{\overline{\gind}}(\gamma_j) \ge \ell_{\overline{\gind}}((\gamma_j)_{[0,S]}) \stackrel{j \to \infty}{\to} \ell_{\overline{\gind}}(\gamma_{|[0,S]}) \stackrel{S \to \infty}{\to} \ell_{\overline{\gind}}(\gamma) = \infty.  
\]
Whence, $T_{\sigma} = \infty$ and thus the Riemannian metric $\overline{\gind}$ is complete. 
 \hfill{$\circledast$}	

Pick a smallest divergent curve $\gamma$ in $(M^n, \overline{\gind})$ issuing from a fixed origin.
Since $(M^n,\gind)$ is complete, reparametrizing $\gamma$ by $\gind$-arclength $s$, it turns out that $\gamma$ is defined for $s \in [0,\infty)$. Whence, because of Claim 1 and since $\gamma$ is $\overline{\gind}$-minimizing between any pair of its points, along $\gamma$ it holds
\[
c_0\int_0^{a} u^{\beta}\psi^2 \di s \le -2t_0\int_{0}^{a}u^\beta \psi\psi_{ss}\di s, \qquad \forall \, a>0
\quad\text{and}\quad \ \psi \in C^2_0([0,a]).
\]
Choosing as test function
$$
\psi(s) = \sin \Big(\frac{\pi s}{a}\Big) \in C^2_0([0,a]),
$$
we get 
\[
\Big(c_0-\frac{2t_0\pi^2}{a^2}\Big) \int_0^a \sin^2 \Big(\frac{\pi s}{a}\Big) u^\beta(\gamma(s)) \di s \le 0,
\]
which gives a contradiction if
$$a > \pi\sqrt{ \frac{2t_0}{c_0}}.$$
This completes the proof of Claim 2.
\end{proof}

\noindent \textbf{Claim 3:} {\em The Riemannian manifold $(M^n,\gind)$ is a totally umbilic sphere.}

\begin{proof}[Proof of Claim 3] 
Since $M^n$ is compact and $\Delta |\Phi|^2 \ge 0$, we deduce that $|\Phi|^2$ is constant. Inequality $u>0$ implies $|\Phi|^2 < b^2$. Thus again from \eqref{lap1} we get 
\[
(b^2-|\Phi|^2)|\Phi|^2 \le 0, 
\]
whence $|\Phi|^2\equiv 0$ and $M^n$ is a sphere.
\end{proof}

This completes the proof of the theorem.

\begin{remark}\label{rem_important}
We point out that Claim 1 can be split into the following two steps:
\begin{itemize}
\item[$i)$] the algebraic estimate 
$$
\Ric+2\beta\abs{\Phi}^2\gind\ge c_0\gind\qquad\text{ for some }c_0\in\R^+,
$$
which in our setting holds for every choice of $\beta>0$. This allows to deduce inequality \eqref{estimate2} from \eqref{elbert};
\item[$ii)$] the inequality
\begin{equation}\label{eqclaim1}
c_0\int_0^au^\beta\psi^2\,\di s\le-2t_0\int_0^au^\beta\psi\psi_{ss}\,\di s \qquad\forall \, \psi\in C_0^2([0,a]),
\end{equation}
along a minimizing $\overline\gind$-geodesic. We prove it when $\beta$ is defined as in \eqref{def_beta}.
\end{itemize}
\end{remark}

\section{The higher-codimensional case}

Let $f:M^n\rightarrow\S^{n+p}$ be a complete submanifold. We denote by $A$ the vector valued second fundamental form of $f$. Consider a local orthonormal frame $\{e_1,\dots,e_{n+p}\}$, with $\{e_1,\dots,e_n\}$ tangent to $M$, and dual coframe $\set{\theta^1,\dots,\theta^{n+p}}$. Using the index convention 
\[
i,j,\ldots\in\set{1,\dots,n}, \qquad \alpha,\beta,\ldots\in\set{n+1,\dots,n+p}, 
\]
$A$ writes in components as
$$
A=h^\alpha_{ij}\theta^i\otimes\theta^j\otimes e_\alpha.
$$
The mean curvature vector is thus defined as 
$$
\mathbf{H}=\frac1nh^\alpha_{kk}e_\alpha,
$$
and the traceless second fundamental form as
$$
\Phi=A-\mathbf Hg=\pa{h^\alpha_{ij}-\frac1nh^\alpha_{kk}\delta_{ij}}\theta^i\otimes\theta^j\otimes e_\alpha=\Phi^\alpha_{ij}\theta^i\otimes\theta^j\otimes e_\alpha.
$$
We will denote by $H$ the norm of $\mathbf{H}$. Observe that, under the assumption that $M$ has parallel mean curvature, namely, that $\nabla^\perp\mathbf{H}=0$, then $H$ turns out to be constant. If $\eta$ is a normal vector field, we denote by $\Phi_\eta$ the bilinear form
$$
\Phi_\eta=\langle\Phi,\eta\rangle=\Phi^\alpha_{ij}\eta^\alpha\theta^i\otimes\theta^j.
$$
A submanifold $M$ will be called pseudo-umbilical if $\mathbf H\neq0$ and $\Phi_{\mathbf H}=0$, i.e. the mean curvature vector lies in an umbilical direction.

We are ready to prove our extension of the results by Simons and Chern, do Carmo \& Kobayashi.

\begin{proof}[Proof of Theorem \ref{teo_minalta}] 
As in \cite{chern}, we can estimate
$$
\frac12\Delta|A|^2\ge\sum_\alpha|\nabla A_{e_\alpha}|^2-|A|^2\set{\pa{2-\frac1p}|A|^2-n}.
$$
Thus, the non-negative function
\[
u= \frac{np}{2p-1} - |A|^2 \qquad \text{satisfies} \qquad \Delta u \le - 2 \frac{2p-1}{p}|A|^2 u \le -2|A|^2 u \quad \text{on } \, M.
\]
If $u(x_0) =0$ for some $x_0 \in M$, by the strong maximum principle $u\equiv 0$, whence
$$
\abs{A}^2\equiv \frac{np}{2p-1}
$$
and the claim follows from \cite{chern}. Although such result is local, Clifford tori and the Veronese surface are connected components of the zero set of some polynomials restricted to the ambient sphere. Hence, the global result follows from the same analyticity argument used in the previous section.

If instead $u>0$ on $M$, then we set $\overline \gind=u^{ 2\beta} \gind$ for some suitable constant $\beta$. To show that $M$ is actually compact, recalling Remark \ref{rem_important} it is enough to show the following \\[0.2cm]
\noindent \textbf{Claim:} {\em For each $\beta>0$, it holds
\[
\Ric +2\beta\abs{A}^2 \gind \ge c_0 \gind
\]
for some $c_0= c_0(\beta,n,p)>0$.} \\[0.2cm]
\emph{Proof of the claim.} Let $X$ be a unit vector, and choose the frame so that $e_1=X$. From the Gauss equation and minimality, 
\[
R_{ik} = (n-1)\delta_{ik}-h^\alpha_{ij}h^\alpha_{jk},
\]
thus Lemma \ref{lem_kato} implies
\begin{equation}\label{riccigammam}
R_{11}=n-1-\sum_\alpha\sum_{j=1}^n\pa{h^\alpha_{1j}}^2 \ge n-1- \frac{n-1}{n}|A|^2.
\end{equation}
Using \eqref{upperA},
$$
n-1\ge\frac{n-1}n\pa{2-\frac1p}\abs{A}^2,
$$
hence for $\tau\in(0,1]$ we have
$$
R_{11}\ge \tau(n-1)+(1-\tau)\frac{n-1}n\pa{2-\frac1p}\abs{A}^2- \frac{n-1}{n}|A|^2,
$$
which gives 
$$
\begin{array}{lcl}
R_{11}+2\beta\abs{A}^2 & \ge & \disp \tau(n-1)+\set{2\beta+\frac{n-1}n \left[(1-\tau)\pa{2-\frac1p}-1\right]} \abs{A}^2 \\[0.5cm]
& \ge & \tau(n-1),
\end{array}
$$
where the last inequality holds for small enough $\tau>0$ depending on $\beta,n$ and $p$.  
\end{proof}

Next, we consider submanifolds with codimension $p \ge 2$ and non-vanishing, parallel mean curvature, where  some subtle difficulties arise. In \cite{santos}, Santos treated the problem for compact submanifolds, obtaining an optimal pinching theorem under the condition that the umbilicity tensor satisfies
\begin{equation}\label{eq_walcy}
\pa{2-\dfrac1{p-1}}|\Phi|^2+\dfrac{n(n-2)}{\sqrt{n(n-1)}}|\Phi_{\bf H}| -n(1+H^2) \le 0.
\end{equation}
Apparently, this is \emph{not} a condition of the type $|\Phi|^2 \le b^2$, so the construction of a conformal factor $u$, if possible, is not evident. To remedy, in the next theorem we slightly strengthen \eqref{eq_walcy} by replacing it with conditions \eqref{condteta} and \eqref{upperphi}. However, in this case we are able to reach the desired conclusions, which rephrase those in \cite{santos}, only in dimension $n \le 6$. 
 
\begin{theorem}\label{teo_cmcalta}
Let $p\ge2$ and let $f:M^n\rightarrow\S^{n+p}$ be a complete, immersed submanifold of dimension $n \le 6$ with parallel, non-zero mean curvature. Assume that the norm of the umbilicity tensor $\Phi$ of $M$ satisfies 
\begin{equation}\label{condteta}
\abs{\Phi_{\mathbf H}}\le\theta H\abs{\Phi}
\end{equation}
for some constant $\theta\in[0,1]$, and
\begin{equation}\label{upperphi}
\abs{\Phi}^2\le b^2,
\end{equation}
where $b=b(n,p,H,\theta)$ is the positive root of the polynomial
$$
P_{n,p,H,\theta}(x)=\pa{2-\dfrac1{p-1}}x^2+\dfrac{n(n-2)}{\sqrt{n(n-1)}}\theta Hx-n(1+H^2).
$$
Then, either $\abs{\Phi}\equiv0$ and $M$ is a totally umbilic sphere, or $\abs\Phi\equiv b$. In this latter case, one of the following occurs:
\begin{itemize}
\item[$i)$] $\theta\in(0,1)$, $p=2$ and $M^n$ covers a $(\theta H)$-torus
 $$
\S^{n-1}(r_1)\times\S^1\pa{r_2}\subset\S^{n+1}_{1+(1-\theta^2)H^2}\subset\S^{n+2}
$$
with $r_1$, $r_2$ uniquely determined by
$$
\frac{(n-1)r_2^2-r_1^2}{nr_1r_2}\sqrt{1+(1-\theta^2)H^2}=\theta H,
$$
$$
r_1^2+r_2^2=\pa{1+(1-\theta^2)H^2}^{-1};
$$
\item[$ii)$] $\theta=0$, $p=2$, $M^n$ is pseudo-umbilical and covers a minimal Clifford torus in a hypersphere 
$$
\S^k\pa{\sqrt{\frac k{n(1+H^2)}}}\times\S^{n-k}\pa{\sqrt{\frac{n-k}{n(1+H^2)}}}\subset\S_{1+H^2}^{n+1}\subset\S^{n+2},
$$
for some $k \in \{1,\ldots,n-1\}$;
\item[$iii)$] $\theta=0$, $n=2$, $p=3$, $M^2$ is pseudo-umbilical and covers a Veronese surface in a hypersphere $\S^4_{1+H^2}\subset\S^5$.
\end{itemize}
\end{theorem}
\begin{proof}
We can choose a local orthonormal frame in such a way that $e_{n+1}=H^{-1}\mathbf H$. Following the computations in \cite[p. 407-410]{santos}, we have
$$
\begin{array}{lcl}
\dfrac12\Delta\abs{\Phi}^2&\ge&\displaystyle\sum_\alpha\abs{\nabla\Phi_{e_\alpha}}^2-\abs{\Phi}^2\set{\pa{2-\dfrac1{p-1}}\abs{\Phi}^2+\dfrac{n(n-2)}{\sqrt{n(n-1)}}\abs{\Phi_{{\bf H}}}-n(1+H^2)\abs{\Phi}^2}\\[.5cm]
&&+\pa{1-\dfrac{1}{p-1}}\abs{\Phi_{e_{n+1}}}^2\pa{2\abs{\Phi}^2-\abs{\Phi_{e_{n+1}}}^2}.
\end{array}
$$

Using \eqref{condteta} and the fact that
$$
\pa{1-\dfrac{1}{p-1}}\abs{\Phi_{e_{n+1}}}^2\pa{2\abs{\Phi}^2-\abs{\Phi_{e_{n+1}}}^2}\ge0
$$
we get
\begin{equation}
\Delta\abs{\Phi}^2\ge-2\abs{\Phi}^2\set{\pa{2-\dfrac1{p-1}}\abs{\Phi}^2+\dfrac{n(n-2)}{\sqrt{n(n-1)}}\theta H\abs{\Phi}-n(1+H^2)}.
\end{equation}
Let $b$ be the positive root of $P_{n,p,H,\theta}$, namely
$$
b=\frac{p-1}{2p-3}\pa{\sqrt{\frac{n(n-2)^2}{4(n-1)}\theta^2H^2+ \frac{(2p-3)n}{p-1} \pa{1+H^2}}-\frac{n(n-2)}{2\sqrt{n(n-1)}}\theta H}.
$$
Then, under the assumption $|\Phi |^ 2\le b^2$, reasoning as in \eqref{pol_roots} we get
\begin{equation}\label{lap-phi}
\Delta\abs{\Phi}^2\ge-2\abs{\Phi}^2P_{n,p,H,\theta}(\abs\Phi)\ge2\pa{2-\dfrac1{p-1}}\abs{\Phi}^2\pa{b^2-\abs{\Phi}^2}\ge0.
\end{equation}
As in the previous section, let us define the function
\[
u\doteq b^2 - |\Phi|^2 
\]
and observe that it satisfies
\[
u\ge 0 \quad\text{and}\quad \Delta u \le -2 |\Phi|^2 u\qquad \text{on } \, M.
\]
If $u(x_0) =0$ for some $x_0 \in M$, by the strong maximum principle $u\equiv 0$, whence $|\Phi |^2\equiv b^2.$ This implies that all inequalities involved in obtaining \eqref{lap-phi} are actually equalities in this case. In particular, $\abs{\Phi_{e_{n+1}}}=\theta \abs\Phi$. Moreover,
$$
\pa{2- \frac{1}{p-1}} \abs{\Phi}^2=n(1+H^2)-\dfrac{n(n-2)}{\sqrt{n(n-1)}} H\abs{\Phi_{e_{n+1}}}
$$
and
$$
\pa{1-\dfrac{1}{p-1}}\abs{\Phi_{e_{n+1}}}^2\pa{2\abs{\Phi}^2-\abs{\Phi_{e_{n+1}}}^2}=0,
$$
which implies that either $p=2$ or $\abs{\Phi_{e_{n+1}}} \equiv 0$ (equivalently, $\theta=0$). Let us consider the case $\abs{\Phi_{e_{n+1}}} \equiv 0$ first. In this case $M$ is pseudo-umbilical and \eqref{upperphi} reduces to
$$
 \abs{\Phi}^2\le\frac{n(1+H^2)}{2-\frac1{p-1}}.
$$
The claim now follows in this case from item $(ii)$ of Proposition 3.1 in \cite{santos}.\\
Let us now assume $p=2$. In this case, since $e_{n+1}$ is parallel so is $e_{n+2}$, hence the normal bundle has zero curvature. As in the proof of Proposition 3.3 of \cite{santos}, we can therefore find parallel normal vector fields $\xi_1$ and $\xi_2$ such that $\xi_2$ is an umbilic direction, i.e. $\Phi_{\xi_2}=0$. The immersion $f$ can thus be split into the composition $f=g_1\circ g_2$, with $g_1:\S^{n+1}_c\to\S^{n+2}$ totally umbilic and $g_2:M^n\to\S^{n+1}_c$. We have that 
$$
\abs{\Phi}^2=\abs{\Phi_{\xi_1}}^2+\abs{\Phi_{\xi_2}}^2=\abs{\Phi_{\xi_1}}^2.
$$
Moreover, setting $H_i=\langle\mathbf H,\xi_i\rangle$, we have that $\abs{\Phi_{\mathbf H}}=\abs{H_1}\abs \Phi$, so 
\[
\abs{H_1}=\theta H \qquad \text{and} \qquad  H^2=H_1^2+H_2^2. 
\]
By the Gauss equation applied to the immersion $g_1$, we find that $c=1+H_2^2$, so $\abs{\Phi_{\xi_1}}$ satisfies 
$$
\abs{\Phi_{\xi_1}}^2+\dfrac{n(n-2)}{\sqrt{n(n-1)}} \theta H\abs{\Phi_{\xi_1}}-n(c+\theta^2H^2)=0.
$$
Item $(ii)$ of Theorem 1.5 of \cite{alencar} can be applied: $M$ is thus locally a $(\theta H)$-torus in $\S^{n+1}_{1+H_2^2}$, with $H_2=(1-\theta^2)H^2$. 

Again, we point out that such rigidity results are based on \cite{chern, alencar} and are therefore local. Nevertheless, Clifford tori, $H$-tori and the Veronese surface are connected components of the zero set of some polynomials restricted to the sphere, so the global result follows as in the previous sections.

If instead $u>0$ on $M$, we show that $(M^n, \overline{\gind})$ with $\overline \gind=u^{2\beta} \gind$ is compact for some suitable constant $\beta$. We prove \\[0.2cm]
\noindent \textbf{Claim:} {\em For each $\beta \ge n/8$, it holds
\begin{equation}\label{eq_low_ricc}
\Ric +2\beta\abs{\Phi}^2 \gind \ge (n-1)\gind.
\end{equation}
}\begin{proof}[Proof of the claim] Having fixed a unit vector $X$, choose the frame so that $e_1=X$. The Gauss equation implies
\[
R_{ij}=(n-1)\delta_{ij}+h^\alpha_{kk}\pa{\Phi^\alpha_{ij}+\frac1nh^\alpha_{ll}\delta_{ij}}-\pa{\Phi^\alpha_{ik}+\frac1nh^\alpha_{ll}\delta_{ik}}\pa{\Phi^\alpha_{kj}+\frac1nh^\alpha_{ll}\delta_{kj}},
\]
%
thus
\begin{equation}\label{riccigamma}
R_{11}=n-1+\pa{n-2}\frac1nh^\alpha_{kk}\Phi^\alpha_{11} +(n-1)H^2-\sum_\alpha\sum_{j=1}^n\pa{\Phi^\alpha_{1j}}^2.
\end{equation}
For $\eps>0$, Young's inequality allows to write
\begin{equation}\label{eq_young}
\frac1nh^\alpha_{kk}\Phi^\alpha_{11}\ge-\frac{ H^2}{2\eps}-\frac\eps2\sum_\alpha\pa{\Phi^\alpha_{11}}^2.
\end{equation}
Plugging this into \eqref{riccigamma} and using Lemma \ref{lem_kato} we obtain
$$
\begin{array}{lcl}
R_{11} & \ge & \disp n-1+\pa{n-1-\frac{n-2}{2\eps}}H^2-\pa{1+\frac{(n-2)\eps}2}\sum_\alpha\sum_{j=1}^n\pa{\Phi^\alpha_{1j}}^2 \\[0.5cm]
& \ge & \disp n-1+\pa{n-1-\frac{n-2}{2\eps}}H^2 - \pa{1 + \frac{(n-2)\eps}2}\frac{n-1}{n}|\Phi|^2.
\end{array}
$$
Hence, $R_{11} + 2\beta|\Phi|^2 \ge n-1$ follows once we solve
$$
\begin{cases}
\disp n-1 -\frac{n-2}{2\eps}\ge0\\[.5cm]
\disp\frac{2\beta n}{n-1}-1-\frac{(n-2)\eps}2\ge0,
\end{cases}
$$
which amounts to
$$
\frac{n-2}{2(n-1)}\le\eps\le\pa{\frac{2\beta n}{n-1}-1}\frac2{n-2}.
$$
These two conditions are compatible if and only if 
$$
\frac{n-2}{2(n-1)}\le\pa{\frac{2\beta n}{n-1}-1}\frac2{n-2},
$$
which is equivalent to imposing $\beta\ge n/8$, as claimed.
\end{proof}

Next, following Remark \ref{rem_important}, we can couple \eqref{eq_low_ricc} with inequality \eqref{eqclaim1} (that holds for $\beta$ satisfying \eqref{def_beta}) to infer the compactness of $M$ whenever
\[
\frac{n}{8} \le \frac{1}{n-2} \qquad (\text{with $<$ if $n=3$}),
\]
which entails $n \le 4$. To reach the conclusion for each $n \le 6$, we observe that the weight $u^\beta$ in each integral of \eqref{eqclaim1} plays no role in the argument leading to the compactness of $M$, described in Claim 2 in the proof of Theorem \ref{THME}. In other words, one may choose 
$$
\varphi=u^{\frac{\beta+\sigma}2}\psi,\qquad\sigma\in\R
$$
as test function in \eqref{def_varphi} to get an inequality like \eqref{eqclaim1} with $u^\sigma$ in place of $u^\beta$, provided that $\beta$ belongs to a suitable interval $J_\sigma$. Notice that $\gamma$ is still a $\overline\gind$-geodesic with $\overline\gind=u^{2\beta}\gind$. As a matter of fact, the choice $\sigma = 0$ was already considered in \cite{shenye2, catino_ronco}, see also $(4)$ in Remark \ref{rem_mainthm}. In view of \eqref{eqclaim1}, we may apply \cite[Corollary 1]{shenye2} or \cite[Theorem 1.1]{catino_ronco} to conclude that $M$ is compact provided 
\[
\beta < \frac{4}{n-1} \qquad \text{if } \, n \ge 4.
\]
Inequality $\frac{n}{8} < \frac{4}{n-1}$ holds if and only if $n \le 6$, concluding the proof. It turns out that one of the values of $\sigma$ which maximizes $\sup J_\sigma$ is precisely $\sigma=0$. 
\end{proof}

\begin{remark}
The proof of \eqref{eq_low_ricc} differs from the corresponding ones in Theorems \ref{THME} and \ref{teo_minalta}. More precisely, unlike \eqref{eq_key} and due to the presence of the parameter $\theta$, in deriving \eqref{eq_young} we implicitly set $\tau=1$. In dimension $n \ge 7$, this choice is indeed optimal for the range of $\beta$ (and thus leads to no admissible $\beta$ for any $\tau \in (0,1]$) unless $\theta$ is larger than some value depending on $n$, an assumption that we would rather avoid. Note that by setting $\tau = 1$ we make no use of the polynomial $P_{n,p,H,\theta}$. In other words, the inequality $|\Phi|^2 \le b^2$ only appears in the construction of the conformal factor $u$, while it plays no role in getting \eqref{eq_low_ricc}. This was quite unexpected to us. 
\end{remark}

\medskip

\noindent\textbf{Acknowledgements}. 
The authors would like to thank Aldir Brasil for kindly pointing out the references \cite{asperti, alias2}. 


\end{document}